\documentclass[12pt, a4paper,english]{article}
\usepackage[T1]{fontenc}
\usepackage[utf8]{inputenc}
\usepackage{babel}

\usepackage{amsfonts, amsthm, amsmath, amssymb}
\usepackage{textcomp, vmargin, microtype}

\usepackage[colorlinks = true, citecolor = cyan]{hyperref}

\newtheorem{teo}{Theorem}[section]
\newtheorem{theorem}[teo]{Theorem}

\theoremstyle{definition}

\numberwithin{equation}{section}

\newcommand{\dt}[1] {\frac{d^{#1}}{d t^{#1}}} 
\newcommand{\qc}[1] {q^{\scriptscriptstyle(#1)}\displaystyle}
\newcommand{\qce}[1] {q_\lambda^{\scriptscriptstyle(#1)}\displaystyle}
\newcommand{\F}[1] {F^{\scriptscriptstyle(#1)}\displaystyle}

\usepackage{tikz,xcolor,hyperref}
\definecolor{lime}{HTML}{A6CE39}
\DeclareRobustCommand{\orcidicon}{%
	\begin{tikzpicture}
	\draw[lime, fill=lime] (0,0) 
	circle [radius=0.16] 
	node[white] {{\fontfamily{qag}\selectfont \tiny ID}};
	\draw[white, fill=white] (-0.0625,0.095) 
	circle [radius=0.007];
	\end{tikzpicture}
	\hspace{-2mm}
}
\foreach \x in {A, ..., Z}{%
	\expandafter\xdef\csname orcid\x\endcsname{\noexpand\href{https://orcid.org/\csname orcidauthor\x\endcsname}{\noexpand\orcidicon}}
}

\title{Nonlocal constants of motion\\
in Lagrangian Dynamics of any order}

\author{
Gianluca Gorni \orcidB{}\\
Universit\`a di Udine,
Dipartimento di\\ Scienze Matematiche, Informatiche e Fisiche\\
via delle Scienze~208, 33100 Udine, Italy\\
\tt{gianluca.gorni@uniud.it}
\and 
Mattia Scomparin \orcidA{}\\ 
Mogliano Veneto 31021, Italy\\
https://orcid.org/0000-0002-2795-5929\\
\tt{mattia.scompa@gmail.com}
\and
Gaetano  Zampieri \orcidC{}\\
Universit\`a di Verona,
Dipartimento di Informatica\\
strada Le Grazie 15, 37134 Verona, Italy\\
\tt{gaetano.zampieri@univr.it}}

\date{}
\begin{document}
\maketitle

\begin{abstract}
We describe a recipe to generate ``nonlocal'' constants of motion for ODE Lagrangian systems.
As a sample application, we recall a nonlocal constant of motion for dissipative mechanical systems, from which we can deduce global existence and estimates of solutions under fairly general assumptions. 
Then we review a generalization to Euler-Lagrange ODEs of order higher than two, leading to first integrals for the Pais-Uhlenbeck oscillator and other systems.
Future developments may include adaptations of the theory to Euler-Lagrange PDEs.
\end{abstract}


\begin{center}
Dedicated to prof.\ Chaudry Masood Khalique\\ on the occasion of his retirement
\end{center}

\section{Introduction}

We are interested in constants of motion for the \emph{Euler-Lagrange equation}
\begin{equation}
 \frac{d}{dt}\partial_{\dot q}L\bigl(t,q(t),\dot q(t)\bigr)-
   \partial_qL\bigl(t,q(t),\dot q(t)\bigr)=0,
\end{equation}
where $L(t,q,\dot q)$ is a smooth scalar valued Lagrangian function,  $t\in  \mathbb{R}$, $q,\dot q\in\mathbb{R}^n$. In the paper \cite{GZNoether} the first and the last author revisited Noether's Theorem, which links first integrals with symmetries of the Lagrangian~$L$. Leaving aside asynchronous perturbations and boundary terms and other issues, here we single out the following simple result. (Notation: the central dot is the scalar product in $\mathbb{R}^n$).

\begin{theorem}\label{nonlocaltheorem2nd} Let $q(t)$ be a solution to the Euler-Lagrange equation
and let $q_\lambda(t)$, $\lambda\in \mathbb{R}$, be a smooth family of perturbed  motions, such that $q_0(t)\equiv q(t)$.
Then the following function of $t$ is constant
\begin{equation}\label{constantofmotionsecondorder}
 \partial_{\dot q}
  L\bigl(t,q(t),\dot q(t)\bigr)\cdot
  \partial_\lambda q_\lambda(t)
  \big|_{\lambda=0} -
  \int_{t_0}^t  \frac{\partial}{\partial\lambda}
  L\bigl(s,q_\lambda(s),\dot q_\lambda(s)\bigr)
  \big|_{\lambda=0}ds\,.
\end{equation}
\end{theorem}

The proof is straightforward: we just take the derivative of the function in \eqref{constantofmotionsecondorder} and use the Euler-Lagrange equation and reverse the order of a double derivative.

We call \eqref{constantofmotionsecondorder} the \emph{constant of motion associated to the family} $q_\lambda(t)$. For a random family, we may expect the constant of motion to be trivial or inconsequential. In general it is \emph{nonlocal}, which means that its value at a time $t$ depends not only on the current state $(t,q(t),\dot q(t))$ at time~$t$, but also on the whole history between $t_0$ and~$t$.

In the original spirit of Noether's theorem we can concentrate the attention to families $q_\lambda(t)$ which make the integrand in \eqref{constantofmotionsecondorder} vanish whenever $L$ enjoys an invariance property. For instance, for the Lagrangian of a particle in the plane under a central force field 
\begin{equation}
  L(t,q,\dot q):=
  \frac{1}{2}m \lVert\dot q\rVert^2
  -U\bigl(t,\lVert q\rVert\bigr),\quad
  q=(q_1,q_2)\in\mathbb{R}^2,
\end{equation} 
the rotation family
\begin{equation}\label{rotationfamily}
  q_\lambda(t):=\begin{pmatrix}
  \cos\lambda & -\sin\lambda\\
  \sin\lambda & \cos\lambda
  \end{pmatrix} \begin{pmatrix}q_1(t)\\ q_2(t)\end{pmatrix}
\end{equation}
exploits the invariance of $L$ under rotations when plugged into~\eqref{constantofmotionsecondorder}, and leads to the conservation of the angular momentum $\det(q(t),\dot q(t))$, which is definitely a \emph{local} first integral.

In a series of papers we have used Theorem~\ref{nonlocaltheorem2nd} to find numerous \emph{nonlocal} constants of motion which are useful. Here is a partial list.
\begin{itemize}
\item
Consider a particle moving in a time-independent potential field $U(q), q\in \mathbb{R}^n$ and \emph{viscous}, i.e., linear, fluid resistance. We see at once that we have global existence in the future of the solutions to the Euler-Lagrange equation. What about the existence in the past? For a quite natural choice of the family $q_\lambda(t)$, after an integration by parts the integrand in 
\eqref{constantofmotionsecondorder} becomes negative, provided we assume $U\ge 0$. Under these conditions the first term in \eqref{constantofmotionsecondorder} is an increasing function of time, which permits to prove the \emph{global existence in the past}. The paper \cite{GZviscousdissipation} also obtains other estimates for the solutions for this system and for some Lane-Emden equation for which  global existence is proved too. We have picked this example for a detailed exposition in Section~\ref{section:dissipative} below, as an illustration of one way of using our nonlocal constants of motion.

\item General time-independent homogeneous potentials of degree~$-2$, for instance Calogero's potential for the one dimensional $n$-body problem with  inversely quadratic pair potentials,  taken from~\cite{GZNoether}, for which an explicit formula of the time-dependence of the distance from the origin is proved, even
though we do not know the shape of the orbit.

\item The Maxwell-Bloch equations for conservative laser dynamics  taken from~\cite{GZMBconservative}, where \eqref{constantofmotionsecondorder} by derivation permits to separate the equations (in a nonstandard way) into a system exhibiting
``fish dynamics'' and a system with central force.
\end{itemize}

Theorem~\ref{nonlocaltheorem2nd} can be generalized to  \emph{nonvariational} Lagrange equations, as we do in~\cite{GZKilling}, among other results on Killing-type equations, and applications for:
\begin{itemize}
\item A particle under 
a time-independent potential field $U(q), q\in\mathbb{R}^n$, and \emph{hydraulic}, i.e. quadratic, fluid resistance, a result taken from Gorni-Zampieri \cite{GZhydraulicdissipative}, for which the result in dynamics is: if $0\le U\le U_{\sup}<+\infty$, all the solutions  for which the initial kinetic energy is strictly greater than $U_{\sup}$ explode in the past in finite time.
\item The dissipative Maxwell-Bloch equations for laser dynamics, taken from Gorni-Residori-Zampieri \cite{GRZlaserdiss}, for which  a quite natural family $q_\lambda(t)$ yields a constant of motion 
\eqref{constantofmotionsecondorder} which, for a special choice of the parameters, turns out to be a 
genuine first integral $N(t,q,\dot q)$ since the integral term vanishes. The first integral permits some kind of separation of  variables.
\end{itemize}

The paper \cite{GZhydraulicdissipative} presents the result on hydraulic fluid resistance and gives a survey on all  other applications we mentioned till now. The novelty in the present survey is the introduction of Scomparin's generalization of Theorem~\ref{nonlocaltheorem2nd} to higher-order Lagrangian systems, which recently appeared in~\cite{Scomparin}.

For $N=1,2,\ldots$ consider the  \emph{higher-order Euler-Lagrange equation} 
\begin{equation}\label{eq:eom}
  \sum_{k=0}^N (-1)^k \dt{k}\,
  \partial_{\qc{k}} L(t,q,\ldots,\qc{N})=0,
\end{equation}
where the $N^{\text{\tiny th}}$-order Lagrangian $L$ is a smooth function with $t\in \mathbb{R}$, and 
$q,\ldots,\qc{N}\in\mathbb{R}^n$. We use $\qc{k}\equiv d^kq/dt^k$.

Within the higherer-order framework, 
a \textit{first integral} is a smooth function 
\begin{equation*}
K(t,q,\qc{1},\qc{2},\ldots)
\end{equation*}
that is constant along all solutions of the Euler-Lagrange equation \eqref{eq:eom}. The celebrated \textit{Noether's Theorem} establishes a relation between invariance proprieties of a Lagrangian and its first integrals \cite{PhysRev.186.1299,Constantelos:1974qu}.

Scomparin's paper\cite{Scomparin} generalizes Theorem~\ref{nonlocaltheorem2nd} to $N^{\scriptscriptstyle \mathrm{th}}$-order Lagrangians:

\begin{teo}\label{teo:nonloc}
Let $t\mapsto q(t)$ be a solution to the Euler-Lagrange equation for smooth $L(t,q,\ldots,\qc{N})$, and let $q_\lambda(t)$, $\lambda \in \mathbb{R}$, be a smooth family of perturbed motions, such that $q_0(t)\equiv q(t)$. Then the following function of t is constant:
\begin{multline}\label{eq:const}
  \sum_{j=1}^N \sum_{k=0}^{j-1} (-1)^k 
  \frac{d^k}{dt^k}\,
  \partial_{\qc{j}} L(t, q,\ldots,\qc{N})
  \cdot\partial_\lambda \qce{j-k-1}
  \big\rvert_{\lambda=0}-\\
  -\int_{t_0}^t\!
  \frac{\partial}{\partial \lambda} \!L(s,q_\lambda, \ldots, \qce{N})
  \big\rvert_{\lambda=0}ds.
\end{multline}
\end{teo}

A
basic higher-order mechanical system is the \emph{Pais-Uhlenbeck oscillator} \cite{PhysRev.79.145}, whose Lagrangian function is
\begin{equation}\label{LPU}
  L^{\scriptscriptstyle\mathrm{PU}}=
  \tfrac{1}{2}\qc{2}^2-\tfrac{1}{2}(w_1^2+w_2^2)\,
  \qc{1}^2+\tfrac{1}{2}w_1^2w_2^2q^2,
\end{equation}
with  $w_1,w_2>0$.  Its Euler-Lagrange  equation is
\begin{equation}
\qc{4}+ (w_1^2+w_2^2)\qc{2}+w_1^2w_2^2\,q=0.
\end{equation}

Higher-order Lagrangians provide a very large class of models for modified gravity theories \cite{Langlois:2015cwa}, 
quantum-loop cosmologies \cite{Liu:2017puc}, 
and string theories \cite{Simon:1990ic}. 
Approaching higher-order mechanics from a new nonlocal point 
of view provides new perspectives to identify novel 
first integrals without necessarily requiring invariance 
proprieties on the already difficult to investigate structure of higher-order Lagrangians. 
Hopefully Theorem~\ref{teo:nonloc} will provide a valuable tool to give a novel insight into stability proprieties of higher-order models and boundedness of related solutions as is done for second order equations  in Kaparulin~\cite{Kaparulin2015EnergyAS, Kaparulin:2019njc} and in the papers mentioned above. 

Section~\ref{section:firstintegrals} below summarizes the main results in Scomparin's paper \cite{Scomparin}.
 
\section{An application to dissipative dynamics}
\label{section:dissipative}

The results in this section are taken from Gorni-Zampieri \cite{GZviscousdissipation}.

Let us consider the first integral of \emph{energy} which is generally derived in Noether's framework by the use of asynchronous perturbations but can also be treated by means of Theorem~\ref{nonlocaltheorem2nd}. 
For a time independent Lagrangian function $L(t,q,\dot q)= \mathcal{L} (q,\dot q)$,  and the \emph{time-shift} family $q_{\lambda}(t)= q(t+\lambda)$ we have 
\begin{equation}
  \partial_{\lambda}L\bigl(t,q_{\lambda}(t),
  \dot q_{\lambda}(t)\bigr)
  \big|_{\lambda=0}=\partial_{q}\mathcal{L}\cdot \dot q(t)+
  \partial_{\dot q}\mathcal{L}\cdot \ddot q(t)=
  \frac{d}{dt}\mathcal{L}(q(t),\dot q(t)).
\end{equation}
Thus the constant of motion~\eqref{constantofmotionsecondorder} is
\begin{multline*}
  \partial_{\dot q}\mathcal{L}\cdot \dot q(t)-
  \int_{t_0}^t \frac{d}{ds}\mathcal{L}(q(s),\dot q(s))ds  
  =\\
  =\partial_{\dot q}\mathcal{L}\bigl(q(t),\dot q(t)\bigr)\cdot
  \dot q(t)-\mathcal{L}\bigr(q(t),\dot q(t)\bigl)+
  \mathcal{L}\bigl(q(t_0),\dot q(t_0)\bigr)
\end{multline*}
which is energy 
\begin{equation}\label{energy}
  E(q,\dot q)=
  \partial_{\dot q}\mathcal{L}(q,\dot q)\cdot \dot   
  q-\mathcal{L}(q,\dot q)
\end{equation}
up to the trivial additive constand $\mathcal{L}(q(t_0),\dot q(t_0))$.

Now, we turn to dissipation. Consider $k>0$, a smooth potential function $U:\mathbb{R}^n\to\mathbb{R}$ defined on the whole
space. The equation of motion of a particle of mass $m>0$ under this potential and viscous dissipation is 
\begin{equation}\label{viscousdissipativeODE}
m\ddot q=-k\dot q-\nabla U(q),\quad q\in\mathbb{R}^n.
\end{equation}

The energy first integral~\eqref{energy} for $k=0$ is 
\begin{equation}
  E(q,\dot q)=
  \textstyle{\frac12}m \lVert\dot q\rVert^2+U(q).
\end{equation}
For $k>0$ this function decreases along solutions
\begin{equation*}
  \dot E=m \dot q\cdot \frac1m
  \bigl(-k\dot q-\nabla U(q)\bigr)+\nabla U(q)\cdot \dot q
  =-k\lVert\dot q\rVert^2\le 0.
\end{equation*} 

In the sequel we assume that \emph{the potential  is bounded from below}, say $U\ge 0$. Then, for any solution $q(t)$ the velocity    $\dot q(t)$ is bounded in the future:
\begin{equation*}
  \frac{1}{2}m \bigl\|\dot q(t)\bigr\|^2\le 
  \frac{1}{2}m \bigl\|\dot q(t)\bigr\|^2+U\bigl(q(t)\bigr)
  \le\frac{1}{2}m \bigl\|\dot q(t_0)\bigr\|^2+U\bigl(q(t_0)\bigr),
  \quad t\ge t_0,
\end{equation*}
so $q(t)$ is bounded for bounded intervals of time  and we get 
\emph{global existence in the future}. What about the past?
   
Notice that for $k=0$, with no dissipation, we have global existence since the above argument holds in the past too. 

Our \eqref{viscousdissipativeODE} can be seen as the
Euler-Lagrange equation for the Lagrangian 
\begin{equation}\label{Caldirola}
  L=e^{kt/m}\Bigl(\frac{1}{2}m \lVert\dot q\rVert^2
  -U(q)\Bigr).
\end{equation}

It is quite natural to consider the family $q_\lambda(t):= q(t+\lambda\,e^{at})$ with $a\in\mathbb{R}$ new parameter and $q(t)$ solution. Indeed for $a=0$ the family reduces to the time-shift used for energy conservation as $k=0$ and the exponential function is easily inspired by the one in~\eqref{Caldirola}.

Then
\begin{multline*}
  \frac{\partial}{\partial\lambda}
  L\bigl(t,q_\lambda(t),\dot q_\lambda(t)\bigr)
  \Big|_{\lambda=0}=\\
  =\frac{d}{dt}\left(-2e^{(a+\frac{k}{m})t}U
  \bigl(q(t)\bigr)\right)+
  e^{(a+\frac{k}{m})t}
  \Bigl(\bigl(a-\textstyle{\frac{k}{m}\bigr)m\bigl\|
  \dot q(t)\bigr\|^2+2\bigl(a+\frac{k}{m}\bigr)}U
  \bigl(q(t)\bigr)\Bigr)
\end{multline*}
where we eliminated $\ddot q(t)$ using the differential equation. For $a=k/m$ it simplifies and we have the constant of motion 
\begin{multline}\label{constant}
  t\mapsto
  \partial_{\dot q}
  L\bigl(t,q(t),\dot q(t)\bigr)\cdot
  \partial_\lambda q_\lambda(t)
  \big|_{\lambda=0}-
  \int_{t_0}^t
  \frac{\partial}{\partial\lambda}
  L\bigl(s,q_\lambda(s),\dot q_\lambda(s)\bigr)
  \Big|_{\lambda=0}ds=\\
  =e^{2kt/m}\Bigl(m\lVert\dot q(t)\rVert^2
  +2U\bigl(q(t)\bigr)\Bigr)
  +4\,\frac{k}{m}\int^{t_0}_t e^{2ks/m}
  U\bigl(q(s)\bigr)ds.
\end{multline}
Since we assumed $U\ge 0$, the last integral  decreases for $t\le t_0$  and the function
\begin{equation}
  t\mapsto e^{2kt/m}
  \Bigl(m\bigl\|\dot q(t)\bigr\|^2+2U(q(t))\Bigr)
\end{equation}
\emph{increases with $t$ for all $t\le t_0$}.

Finally, we have the estimate for $t\le t_0$:
\begin{equation}
  m\bigl\|\dot q(t)\bigr\|^2\le e^{2k(t_0-t)/m}
  \Bigl(m\bigl\|\dot q(t_0)\bigr\|^2+2U\bigl(q(t_0)\bigr)\Bigr).
\end{equation}
In a bounded interval $(t_1,t_0]$ the velocity $\dot q(t)$ is bounded, so also $q(t)$  and we have global existence of solutions. Summing up:

\begin{teo}\label{teo:globalexistence}
If $k>0$ and $U$ is a smooth potential on~$\mathbb{R}^n$ which is bounded from below, then all solutions of the dissipative equation $\ddot q=-k\dot q-\nabla U(q)$ are defined for all~$t\in\mathbb{R}$.
\end{teo}

\section{First integrals for higher-order Lagrangians}\label{section:firstintegrals}

Generally, we cannot expect that Theorem \ref{teo:nonloc} yields true first integrals for a random choice of the family $q_\lambda(t)$. However, few and precious Lagrangians make our machinery work. The simplest example is for autonomous Lagrangians:

\begin{theorem}\label{teo:timeshift}
Let $t\mapsto q(t)$ be a solution to the Euler-Lagrange equation for smooth time-independent 
$\mathcal{L}(q,\ldots,\qc{N})$. Then the following function is a first integral:
\begin{multline}\label{eq:timeshift}
  K_1(q,\ldots,\qc{2N-1})=\\
  =\sum_{i=1}^N
  \sum_{k=0}^{i-1} (-1)^k   \dt{k}
  \partial_{\qc{i}} \mathcal{L}(q,\ldots,\qc{N})
  \cdot\qc{i-k}-\mathcal{L}(q,\ldots,\qc{N})\,.
\end{multline}
\end{theorem}

It is important to notice that Theorem \ref{teo:timeshift} recovers the Noetherian result of \cite{Logan:1975au} for $N=2$ Lagrangians.

Since the Pais-Uhlenbeck Lagrangian $L^{\scriptscriptstyle\mathrm{PU}}$ of formula~\eqref{LPU} is time-independent, we deduce the following first integral:  
\begin{equation}
  2K^{\scriptscriptstyle\mathrm{PU}}_1=
  \qc{2}^2-(w_1^2+w_2^2)\qc{1}^2-
  2\qc{3}\qc{1}-w_1^2w_2^2q^2.
\end{equation}

In \cite{GZNoether} the authors proved energy conservation for the canonical harmonic oscillator starting from nonlocal space-changes. Consequently, using Theorem~\ref{teo:nonloc}, we deduce that first integrals are easy to be found if 
$\partial_{\qc{j}} L\!\propto d^j\partial_q L/dt^j$ for all $j=1,\ldots,N$:

\begin{theorem}\label{teo:nonlocspace}
Consider a smooth Lagrangian $L(t,q,\ldots,\qc{N})$ for which there exists a set of constant parameters $\rho_1\ldots\rho_N\in\mathbb{R}$ such that
\begin{equation}
  \partial_{\qc{i}} L=\rho_i \frac{d^i}{dt^i}
  \partial_q L\qquad\text{for all motions and }
  i\in\{1,\ldots,N\}.
\end{equation} 
Let $t\mapsto q(t)$ be a solution to the Euler-Lagrange equation, and define 
$F^{\scriptscriptstyle(0)}=\sum_{j=1}^N(-1)^{j+1} d^{j-1} \partial_{\partial \qc{j}} L/dt^{j-1}$.
Then 
$F^{\scriptscriptstyle(\ell)}= d^{\ell-1}\partial_qL/dt^{\ell-1}$
 with $\ell \in \{1,\ldots, 2N\}$, and the following function is a first integral:
\begin{multline}\label{eq:rexx}
  K_2(t,q,\ldots,\qc{3N-1})=\\
  =\sum_{i=1}^N\rho_i
  \Biggl[\sum_{k=0}^{i-1}(-1)^k\F{i+k+1}\cdot \F{i-k-1}-
  \frac{1}{2}\lVert\F{i}\rVert^{2}
  \Biggr]-
  \frac{1}{2}\lVert\F{0}\rVert^{2}.
\end{multline}
\end{theorem}

Notice that $L^{\scriptscriptstyle\mathrm{PU}}$ satisfies Theorem \ref{teo:nonlocspace} with $\rho_1=-(w_1^2+w_2^2)w_1^{-2}w_2^{-2}$ and $\rho_2=w_1^{-2}w_2^{-2}$. Hence, $\F{0} =-(w_1^2+w_2^2)\qc{1} - \qc{3}$ and $\F{\ell}=w_1^2w_2^2\qc{\ell-1}$ ($\ell=0,1,2$) give
\begin{multline}
  2K^{\scriptscriptstyle\mathrm{PU}}_2=
  \bigl(w_1^4+w_1^2w_2^2+ w_2^4\bigr)  
  \qc{1}^2+\qc{3}^2+2w_1^2w_2^2\,q\qc{2}+\\
  +(w_1^2+w_2^2) \bigl(2\qc{3}\qc{1}+ w_1^2w_2^2q^2\bigr)
\end{multline}
If combined, $K_1^{\scriptscriptstyle\mathrm{PU}}$ and $K_2^{\scriptscriptstyle\mathrm{PU}}$ recover the two first integrals recently proposed by \cite{Kaparulin:2019njc}.

Another interesting situation generating first integrals arises when the Lagrangian, evaluated on a perturbed motion $q_\lambda(t)$, has constant derivative at $\lambda = 0$. 
\begin{theorem} \label{teo:simm}
Consider a smooth Lagrangian, and suppose that for a given family of perturbed motions $q_\lambda$ there exists a constant  $\mu\in \mathbb{R}$ such that 
$\partial_\lambda L(t,q_\lambda,\ldots,q_{\lambda}^{\scriptscriptstyle{(N)}})\big\rvert_{\lambda=0}\!=\!\mu\,.$
Then, the following function is a first integral:
\begin{multline}\label{eq:lxc}
  K_3(t,q,\ldots,\qc{2N-1})=\\
  =\sum_{i=1}^N \sum_{k=0}^{i-1} (-1)^k \dt{k}\partial_{\qc{i}} 
  L(t,q,\ldots,\qc{N})\cdot
  \partial _\lambda \qce{i-k-1}\big|_{\lambda=0}-\mu t\,.
\end{multline}
\end{theorem}

Let $q=(q_1,q_2)\in\mathbb{R}^2$ and consider the rotation family $q_\lambda(t)$ in \eqref{rotationfamily}. When evaluated on $q_\lambda(t)$, $L^{\scriptscriptstyle\mathrm{PU}}$ does not depend on $\lambda$, hence  Theorem \ref{teo:simm} gives the new angular-momentum-like first integral 
\begin{equation}
  K^{\scriptscriptstyle\mathrm{PU}}_3=
  (w_1^2+w_2^2)\det\bigl(\qc{1},q\bigr)+
  \det\bigl(\qc{1},\qc{2}\bigr)+\det\bigl(\qc{3},q\bigr).
\end{equation}

In Scomparin \cite{Scomparin}, the machinery is also applied to analyze a higher-order generalization of the Pais-Uhlenbeck oscillator \cite{Andrzejewski:2014zia,Boulanger:2018tue}, and  a simple Degenerate Higher-Order model of Scalar-Tensor (DHOST) theory that, in these last years, inspired many modified gravity theories \cite{Langlois:2015cwa}. Again, the full consistency of the machinery is confirmed.

\section*{Acknowledgments}

Paper written under the auspices of INdAM (Istituto Nazionale di Alta Matematica).


\end{document}